\documentclass{amsart}
\usepackage{wasowicz}

\info{%
 Published: J. Appl. Anal. 1~(1995), 173--179\\
 A manuscript was originally written in 1993.\\
 This is a version rewritten in 2008 using the \texttt{amsart} document class.}

\begin{document}
\author{\SW}
\address{\SWaddr}
\email{\SWmail}
\title{On affine selections of set--valued functions}
\date{1993}
\begin{abstract}
 The main result of this paper is the theorem stating that every convex
 set--valued function $F:I\mapsto\compact(Y)$, where $I\subset\R$ is an
 interval and $Y$ is locally convex space, possesses an affine selection.
 In the case if $Y=\R$ and the values of $F$ are closed real intervals
 we can replace the assumption of convexity of $F$ by the more general
 condition.
\end{abstract}
\maketitle

\section{Introduction}

K.~Nikodem and Sz.~Wąsowicz have proved (cf.\ \cite[theorem 1]{NikWas95})
that if the functions $f,g:I\mapsto\R$, where $I\subset\R$ is an interval,
fulfil for every $x,y\in I$, $t\in [0,1]$ the following condition
\[
\begin{aligned}\label{nasz}
 f\bigl(tx+(1-t)y\bigr)&\le tg(x)+(1-t)g(y),\\
 g\bigl(tx+(1-t)y\bigr)&\ge tf(x)+(1-t)f(y),
\end{aligned}
\]
then there exists an affine function $h:I\to\R$ such that $f\le h\le g$ on $I$.
The simple consequence of this fact is Theorem~\ref{tw1} which we prove at the
beginning of this paper. Next we prove that every convex set--valued function
$F:I\to\compact(Y)$, where $Y$ is a~locally convex space and $\compact(Y)$ is the family of
all compact non--empty subsets of~$Y$, possesses an affine selection. In the case
$Y=\R^{n}$, $n\in\N$, we also present a~direct inductive proof of this theorem.

\subsection*{Notation}
By $I$ we will denote any fixed real interval. If $X$ is a~topological vector space
then we admit the following notation:
\begin{align*}
  \nonempty(X)&=\{\,A\subset X\,:\,A\ne\emptyset\,\}\\
   \compact(X)&=\{\,A\in\nonempty(X)\,:\,A\;\text{ is a compact set}\,\}\\
 \compactconvex(X)&=\{\,A\in\compact(X)\,:\,A\;\text{ is a convex set}\,\}.
\end{align*}
The term \emph{set--valued function} will be abbreviated in the form
\emph{s.v.\ function}.

If $X,Y$ are vector spaces and $D\subset X$ is a convex set then we say that
a~s.v. function $F:D\to\nonempty(Y)$ is
\begin{enumerate}[(a)]
 \item
  \emph{convex}, if for every $x,y\in D$, $t\in [0,1]$
  \[
   tF(x)+(1-t)F(y)\subset F\bigl(tx+(1-t)y\bigr);
  \]
 \item
  \emph{concave}, if for every $x,y\in D$, $t\in[0,1]$
  \[
   F\bigl(tx+(1-t)y\bigr)\subset tF(x)+(1-t)F(y).
  \]
\end{enumerate}

A function $f:D\to Y$ is called the \emph{selection} of the~s.v.\ function
$F:D\to\nonempty(Y)$ iff $f(x)\in F(x)$ for every $x\in D$.

\section{Results}

Let us start with the translation of the theorem mentioned in the
Introduction into the~s.v.\ functions language
(cf.~\cite[Remark 14]{ISFE31}).
\begin{thm}\label{tw1}
 A s.v.\ function $F:I\to\compactconvex(\R)$ possesses an affine selection
 iff for every $x,y\in I$, $t\in [0,1]$ the following condition holds
 \begin{equation}\label{nasz_warunek}
  F\bigl(tx+(1-t)y\bigr)\cap [tF(x)+(1-t)F(y)]\ne\emptyset.
 \end{equation}
\end{thm}
Before we start the proof let us observe that \eqref{nasz_warunek} is the
weakest condition guaranteeing the existence of an affine selection for the
s.v.\ function~$F$.
\begin{proof}[Proof of Theorem~\ref{tw1}]
 If there exists an affine selection $f:I\to\R$ of the s.v.\ function $F$ then
 the condition \eqref{nasz_warunek} is obvious.

 Let us assume that the condition \eqref{nasz_warunek} holds for every
 $x,y\in I$, $t\in [0,1]$. Let
 \[
  f(x):=\inf F(x),\quad g(x):=\sup F(x),\qquad x\in I.
 \]
 Then
 \begin{equation}\label{w3}
  F(x)=\bigl[f(x),g(x)\bigr],\quad x\in I
 \end{equation}
 We will show that $f,g:I\to\R$ fulfil the condition \eqref{nasz}
 for any fixed $x,y\in I$, $t\in [0,1]$.
 Let $z\in F\bigl(tx+(1-t)y\bigr)\cap\bigl[tF(x)+(1-t)F(y)\bigr]$.
 There exist $z_{1}\in F(x)$, $z_{2}\in F(y)$ such that $z=tz_{1}+(1-t)z_{2}$.
 Using the definitions of $f$ and $g$ we get
 \begin{align*}
  f\bigl(tx+(1-t)y\bigr)&\le z=tz_{1}+(1-t)z_{2}\le tg(x)+(1-t)g(y),\\
  g\bigl(tx+(1-t)y\bigr)&\ge z=tz_{1}+(1-t)z_{2}\ge tf(x)+(1-t)f(y).
 \end{align*}
 Then there exists an affine function $h:I\to\R$ such that
 \begin{equation}\label{w4}
  f(x)\le h(x)\le g(x),\quad x\in I
 \end{equation}
 (cf.~\cite[Theorem 1]{NikWas95}). Conditions \eqref{w3}, \eqref{w4} imply that
 $h(x)\in F(x)$, $x\in I$, which completes the proof.
\end{proof}
\begin{rem}\label{uw1}
 It is well known that if a s.v.\ function $F:I\to\compactconvex(\R)$ is convex (or concave)
 then $F$ has an affine selection. Applying the above theorem we get the new
 proof of this fact.
\end{rem}
\begin{rem}
 The assumption of compactness of the sets $F(x)$, $x\in I$ in
 Theorem~\ref{tw1} is essential. Consider two s.v.\ functions
 $F:\R\to\nonempty(\R)$ and $G:(-1,1)\to\nonempty(\R)$ defined by the formulas
 \begin{align*}
  F(x)&=[x^{2},+\infty),\quad x\in\R,\\
  G(x)&=(x^{2},1),\quad x\in (-1,1).
 \end{align*}
 It is easy to see that $F$ and $G$ are convex, but they do not have any
 affine selection.
\end{rem}
\begin{rem}\label{Ela}
 It is known (cf.~\cite[Remark 1]{Nik89}, \cite[Remark 2]{NikWas95}) that
 a~s.v.\ function $F:D\to\compactconvex(\R)$, where $D\subset\R^{2}$ is
 a~convex set, need not possess any affine selection although~$F$ fulfils
 \eqref{nasz_warunek}. We can also find the example of the s.v.\ function
 $F:I\to\compactconvex(\R^{2})$ which fulfils \eqref{nasz_warunek} and does not have
 any affine selection. So, Theorem~\ref{tw1} can not be generalized both
 for s.v.\ functions defined on the convex subset of the plane and for
 s.v.\ functions $F:I\to\compactconvex(\R^{2})$. Below we present an example which
 is due to E. Sadowska from Bielsko--Biała.

 Let $I=[0,4]$ (only in this remark) and $F:I\to\compactconvex(\R^{2})$ be
 defined as follows
 \begin{align*}
  F(0)&=[-4,4]\times\{1\},\\
  F(1)&=\{-1\}\times [-4,4],\\
  F(2)&=[-4,4]\times\{-1\},\\
  F(3)&=\{1\}\times [-4,4],\\
  F(4)&=\{\,(x,x)\,:\,x\in [-4,4]\,\},\\
  F(x)&=[-4,4]\times [-4,4]\text{ for all } x\in I\setminus\{0,1,2,3,4\}.
 \end{align*}
 One can prove that $F$ fulfils \eqref{nasz_warunek}. We will show that~$F$
 does not possess any affine selection. An easy computation shows that the
 straight line
 \[
  \ell\,:\,\left\{\begin{aligned}
                  x&=1-\xi\\
                  y&=-1-\xi\\
                  z&=\xi
             \end{aligned}\right.,\quad\xi\in\R
 \]
 is the only line which intersects four segments $F(0)$, $F(1)$, $F(2)$
 and $F(3)$. But $\ell$ does not intersect the segment $F(4)$. So the
 s.v.\ function $F$~has not any affine selection.
\end{rem}

It is well known that every continuous function $f:I\to I$ has a
fixed point (if $I$ is the closed interval). On the second hand, every
affine function $f:I\to I$ is continuous. So, as a~consequence of
Theorem~\ref{tw1} we obtain the following
\begin{cor}
 If the interval $I$ is closed then every s.v.\ function $F:I\to\compactconvex(I)$
 fulfilling \eqref{nasz_warunek} has a~fixed point
 (i.e.\ there exists a~point $x\in I$ such that $x\in F(x)$).
\end{cor}

Now we shall prove the main theorem of this paper. We present two proofs. One
of them is an application of K. Nikodem's results (cf. \cite{Nik89}) and it
requires the Axiom of Choice. The second one is direct and inductive
but it works for finite--dimensional spaces $\R^{n}$.

\begin{thm}\label{tw2}
 Let $Y$ be a locally convex topological vector space. Every convex
 s.v.\ function $F:I\to\compact(Y)$ possesses an affine selection $f:I\to Y$.
\end{thm}

\begin{proof}
 Let $F:I\to\compact(Y)$ be a~convex s.v.\ function. Then, in particular,
 $F$ is a~midconvex s.v.\ function (it fulfils the condition of convexity
 with $t=1/2$), and so $F$ possesses a~Jensen selection $f:I\to Y$
 i.e.\ $f(\frac{x+y}{2})=\frac{f(x)+f(y)}{2}$, $x,y\in I$
 (cf.\ \cite[Lemma 2]{Nik89}).

 Since $F$ is continuous on $\Int I$ as a convex s.v.\ function defined on
 a~subset of $\R$ (cf.\ \cite[Theorem 3.7]{Nik89ZN}), also $f$ is
 continuous (cf.\ \cite[Theorem 4.3]{Nik89ZN} with $K=\{0\}$, $G=F$ and $F=f$).
 Thus $f$ as a~continuous Jensen function is an affine function, which
 completes the proof.
\end{proof} 

 In the above proof we have used K. Nikodem's results which require the
 Lemma of Kuratowski--Zorn and some versions of the Theorem of
 Hahn--Banach. Below we give an inductive proof in the case
 $Y=\R^{n}$, $n\in\N$.
 
\begin{proof}[Second proof (for $Y=\R^{n}$)]
 Before we start an induction on $n$ let us notice that if $F$
 is convex then the values of $F$ are convex sets.
 
 If $n=1$ then our theorem follows directly from Remark~\ref{uw1}.
 Assume that every convex s.v.\ function $G:I\to\compact(\R^{n})$ has an affine
 selection $g:I\to\R^{n}$. Let $F:I\to\compact(\R^{n+1})$ be a convex
 s.v.\ function. For any $x\in I$ we put
 \[
  G(x):=\bigl\{\,y\in\R^{n}\,:\,\exists_{z\in\R}\,(y,z)\in F(x)\,\bigr\}
 \]
 It is easy to verify that $G(x)$ is a compact and non--empty subset of
 $\R^{n}$, i.e.\ $G:I\to\compact(\R^{n})$. We will check that $G$
 is a convex s.v.\ function. Fix any $x_{1},x_{2}\in I$, $t\in [0,1]$. Let
 $y\in tG(x_{1})+(1-t)G(x_{2})$.
 There exist $y_{i}\in G(x_{i})$, $i=1,2$, such that $y=ty_{1}+(1-t)y_{2}$.
 So there exist $z_{i}\in\R$, $i=1,2$, such that
 $(y_{i},z_{i})\in F(x_{i})$, $i=1,2$.
 Let $z=tz_{1}+(1-t)z_{2}$. Since $F$ is convex we get
 \begin{multline*}
  (y,z)=t(y_{1},z_{1})+(1-t)(y_{2},z_{2})\\
  \in tF(x_{1})+(1-t)F(x_{2})\subset F(tx_{1}+(1-t)x_{2}).
 \end{multline*}
 So we obtain $y\in G(tx_{1}+(1-t)x_{2})$.

 Let $g:I\to\R^{n}$ be an affine selection of $G$. Let us define
 \[
  H(x):=\bigl\{\,z\in\R\,:\,\bigl(g(x),z\bigr)\in F(x)\,\bigr\},\quad x\in I.
 \]
 Obviously, $H(x)$ is a compact and convex subset of $\R$. From the fact that
 $g(x)\in G(x)$ we get that $H(x)\ne\emptyset$, $x\in I$. So $H:I\to\compactconvex(\R)$.
 It is not difficult to check that $H$ is a convex s.v.\ function.

 Let $h:I\to\R$ be an affine selection of $H$. By putting
 $f(x):=\bigl(g(x),h(x)\bigr)$, $x\in I$, we obtain an affine selection
 of the s.v.\ function $F$, which completes the proof.
\end{proof} 
\begin{rem}
 The above proof was obtained as a result of investigations if
 s.v.\ functions $F:I\to\compactconvex(\R^{n})$ fulfilling the condition
 \eqref{nasz_warunek} possess an affine selection. First the author
 hoped that the method of decreasing of the dimension will give an effect
 in this case. But the example in Remark~\ref{Ela} gave the negative answer
 of the mentioned problem. However in the case of convex s.v.\ functions this
 direct method gives the simple way of constructing the affine selections.
\end{rem}
\begin{rem}
 It is known that convex s.v.\ functions $F:D\to\compactconvex(\R)$, where
 $D$ is a convex subset of $\R^{n}$, $n\ge 2$, need not possess any
 affine selection (cf.\ \cite[Remark 1]{Nik89}). However, if $D$ is a~convex
 cone with base in a real linear space, then every convex s.v.\ function
 defined on $D$ with compact values in a real locally convex space has an
 affine selection. It was obtained recently by A.~Smajdor and
 W.~Smajdor~\cite{SmaSma96}.
\end{rem}
\begin{rem}
 It is known that the convex s.v. function $F:I\to\compact(Y)$, where $Y$ is
 a~topological vector space are continuous. An application of our results gives
 another proof of this fact in the case when $Y=\R^{n}$.
 Using Theorem~\ref{tw2} (for $Y=\R^{n}$) we obtain the existence of an affine
 selection $f$ of $F$. Since $f$ is affine and $f:I\to\R^{n}$, it must be
 continuous. Then $F$ has a~continuous selection, so $F$ must be continuous
 on $I$ (cf.\ \cite[Theorem 3.3]{Nik89ZN}).
\end{rem}

\bibliographystyle{amsplain}
\bibliography{was_own,was_pub}

\end{document}